\documentclass{article}
\usepackage{amssymb,amsfonts,amsmath,amsthm,enumerate}
\usepackage[numeric]{amsrefs}
\usepackage{mathtools}
\usepackage[nocompress]{cite}

\def \C{{\mathbb C}}

\def \f{\frac}
\def \t{\theta}
\def \z{\zeta}

\def \Q{{\mathbb Q}}

\def \R{{\mathbb R}}
\def \Re{\operatorname{Re}}

\def \Z{{\mathbb Z}}
\def \l{\left}
\def \r{\right}

\def \z{\zeta}

\newtheorem{theorem}{Theorem}
\newtheorem{thm}{Theorem}

\numberwithin{equation}{section} \theoremstyle{definition}

\begin{document}

\title{Rigidity of Euler Products}
\author{Shin-ya Koyama\footnote{Department of Biomedical Engineering, Toyo University,
2100 Kujirai, Kawagoe, Saitama, 350-8585, Japan.} \ \& Nobushige Kurokawa\footnote{Department of Mathematics, Tokyo Institute of Technology, 
Oh-okayama, Meguro-ku, Tokyo, 152-8551, Japan.}}

\maketitle

\begin{abstract}
We report a simple rigidity theorem for certain Euler products.
\end{abstract}



\section*{Introduction}
For purely imaginary numbers $a,$ $b$, $c\in i\R$ we study the meromorphy of the associated Euler product
$$
Z^{abc}(s)=\prod_{p:\text{ prime }}(1-(p^a+p^b)p^{-s}+p^{c-2s})^{-1}
$$
in the family of Euler products
$$
\frak Z^{ab}=\{Z^{abc}(s)\ |\ c\in i\R\}.
$$
This family $\frak Z^{ab}$ contains
$$
\z(s-a)\z(s-b)=Z^{ab(a+b)}(s),
$$
which is a meromorphic function in all $s\in\C$
(with a functional equation under $s\longleftrightarrow 1+a+b-s$).
We prove the converse:
\begin{theorem}
If $Z^{abc}(s)$ is meromorphic in all $s\in \C$, we have $a+b=c$ and 
$$
Z^{abc}(s)=\z(s-a)\z(s-b).
$$
\end{theorem}

This shows rigidity of the Euler product $\z(s-a)\z(s-b)$ in the family
$\frak Z^{ab}$ concerning the meromorphy on the entire $\C$.

The next result gives a detailed meromorphy for $a+b\ne c$.
\begin{theorem}
If $a+b\ne c$, then $Z^{abc}(s)$ has an analytic continuation to $\Re(s)>0$ as a meromorphic function with the natural boundary
$\Re(s)=0$.
More precisely, each point on $\Re(s)=0$ is a limit point of poles of 
$Z^{abc}(s)$ in $\Re(s)>0$.
\end{theorem}

We notice generalizations in \S4 in the text.
Our theorems follow from results of Kurokawa \cite{K1, K2, K3} extending results of Estermann \cite{E}.

We remark that our result characterizes $\z(s-a)\z(s-b)$ by the meromorphy in all $s\in\C$
only in contrast to usual ``converse theorems''
originated by Hamburger \cite{H} and Hecke \cite{He} where the functional equation and the attached automorphic
form are important; $\z(s-a)\z(s-b)$ corresponds to a Maass wave form studied by Maass \cite{M}.

\section{Euler datum}
We use the triple $E=(P,\ \R,\ \alpha)$, where $P$ is the set of all prime numbers,
$\R$ denotes the real numbers, and $\alpha$ is the map
$\alpha:\ P\to\R$ given by $\alpha(p)=\log p$. 
Such a triple is a simple example of Euler datum studied in \cite{K1, K2, K3};
generalized Euler data are treated there with general ``primes'' $P$
and general topological groups $G$ instead of $\R$.

Let $R(\R)$ be the virtual character ring of $\R$ defined as
$$
R(\R)=\l\{\l.\sum_{a\in i\R}m(a)\chi_a\ \r|\ m(a)\in\Z,\ m(a)=0\text{ except for finitely many }a\r\},
$$
where $\chi_a$ is a (continuous) unitary character
$\chi_a:\ \R\to U(1)$ given by $\chi_a(x)=e^{ax}$ for $x\in\R$.

For a polynomial 
$$
H(T)=\sum_{m=0}^n h_m T^m\in1+TR(\R)[T]
$$
we denote by $L(s,E,H)$ the Euler product
$$
L(s,E,H)=\prod_{p\in P}H_{\alpha(p)}(p^{-s})^{-1},
$$
where
$$
H_x(T)=\sum_{m=0}^n h_m(x)T^m\in 1+T\C[T].
$$
For example, let $a$, $b$, $c\in i\R$, then the polynomial
$$
H^{abc}(T)=1-(\chi_a+\chi_b)T+\chi_cT^2
$$
in $1+TR(\R)[T]$ gives the Euler product
$$
L(s,E,H^{abc})=\prod_{p\in P}\l(1-(p^a+p^b)p^{-s}+p^{c-2s}\r)^{-1}
$$
since
\begin{align*}
H^{abc}_{\alpha(p)}
&=1-(\chi_a(\log p)+\chi_b(\log p))p^{-s}+\chi_c(\log p)p^{-2s}\\
&=1-(p^a+p^b)p^{-s}+p^{c-2s}.
\end{align*}

\section{Unitariness and meromorphy}
Let $E=(P,\ \R,\ \alpha)$ as in \S1 and take a polynomial $H(T)$ in $1+TR(\R)[T]$ of degree $n$.
We say that $H(T)$ is {\it unitary} when there exist functions $\theta_j;\ \R\to\R$ satisfying
$$
H_x(T)=(1-e^{i\theta_1(x)}T)\cdots(1-e^{i\theta_n(x)}T)
$$
for all $x$.

The main theorem proved in \cite{K2} gives in this particular situation the following result.
\begin{thm}\mbox{}\\
\begin{enumerate}[\rm(1)]
\item 
\vskip -8mm
If $H(T)$ is unitary, then $L(s,E,H)$ is meromorphic in all $s\in\C$.
\item
If $H(T)$ is not unitary, then $L(s,E,H)$ is meromorphic in $\Re(s)>0$ with the natural boundary.
Moreover, each point on $\Re(s)=0$ is a limit point of poles of $L(s,E,H)$ in $\Re(s)>0$.
\end{enumerate}
\end{thm}
This theorem was proved in \cite{K2} (p.45, \S8, Theorem1) since our
$E=(P,\ \R,\ \alpha)$ is nothing but $\overline{E_0(\Q/\Q)}$ there.

\section{Proof of rigidity}
After looking Theorem 1 recalled in \S2 we see that Theorems A and B in Introduction
are both derived from the following result:
\begin{thm}
Let $a$, $b$, $c\in i\R$ and
$$
H^{abc}(T)=1-(\chi_a+\chi_b)T+\chi_cT^2\in1+TR(\R)[T].
$$
Then the following conditions are equivalent.
\begin{enumerate}[\rm(1)]
\item \vskip -4mm
$a+b=c$
\item $H(T)$ is unitary.
\end{enumerate}
\end{thm}

\noindent
{\it Proof.}
(1)$\Longrightarrow$(2):
From $a+b=c$ we get
\begin{align*}
H^{abc}(T)&=1-(\chi_a+\chi_b)T+\chi_a\chi_b T^2\\
&=(1-\chi_a T)(1-\chi_b T).
\end{align*}
This gives
\begin{align*}
H_x^{abc}(T)&=(1-\chi_a(x) T)(1-\chi_b(x) T)\\
&=(1-e^{ax}T)(1-e^{bx} T)
\end{align*}
for $x\in\R$. Hence $H^{abc}(T)$ is unitary, by
$$
|e^{ax}|=|e^{bx}|=1.
$$

\noindent
(2)$\Longrightarrow$(1):
Assume that $H^{abc}(T)$ is unitary, and set
$$
H_x^{abc}(T)=(1-e^{i\t_1(x)}T)(1-e^{i\t_2(x)} T)
$$
with $\t_j:\ \R\to\R$. Then comparing with
$$
H_x^{abc}(T)=1-(e^{ax}+e^{bx})T+e^{cx} T^2
$$
\makeatletter\tagsleft@true\makeatother
we obtain
\begin{align}
e^{ax}+e^{bx}&=e^{i\t_1(x)}+e^{i\t_2(x)},\tag{3.1}\\
e^{cx}&=e^{i(\t_1(x)+\t_2(x))}.\tag{3.2}
\end{align}
Note that the complex conjugation of (3.1) gives
$$
e^{-ax}+e^{-bx}=e^{-i\t_1(x)}+e^{-i\t_2(x)}.
$$
Since
$$
e^{-ax}+e^{-bx}=e^{-(a+b)x}(e^{ax}+e^{bx})
$$
and
$$
e^{-i\t_1(x)}+e^{-i\t_2(x)}=e^{-i(\t_1(x)+\t_2(x))}(e^{i\t_1(x)}+e^{i\t_2(x)})
$$
we obtain the equality
$$
e^{-(a+b)x}(e^{ax}+e^{bx})=e^{-cx}(e^{ax}+e^{bx})
$$
by using (3.1) and (3.2).

Hence we get
$$
(e^{(a+b-c)x}-1)(e^{ax}+e^{bx})=0
$$
for all $x\in\R$. Especially
$$
\f{e^{(a+b-c)x}-1}x(e^{ax}+e^{bx})=0
$$
for all $x\in\R\setminus\{0\}$. Thus letting $x\to0$ we obtain the desired equality
$$
a+b-c=0.
$$
\hfill\qed

\section{Generalizations}
From the proof above it would be easy to see that we have generalizations of Theorems A and B
by using results of \cite{K1, K2, K3}. Hence we notice simple results only.

\subsection{Dedekind case.}
Let $\z_F(s)$ be the Dedekind zeta function of a finite extension field $F$ of the rational number field $\Q$.
Let $a,$ $b$, $c\in i\R$ and
$$
Z_F^{abc}(s)=\prod_{P\in\mathrm{Specm}(O_F)}
(1-(N(P)^a+N(P)^b)N(P)^{-s}+N(P)^{c-2s})^{-1},
$$
where $P$ runs over the set $\mathrm{Specm}(O_F)$ of maximal ideals of the integer ring $O_F$ of $F$.
Then we have exactly the same Theorems A and B chracterizing 
$\z_F(s-a)\z_F(s-b)$ among $Z_F^{abc}(s)$ by using Theorem 1 of \cite[\S8]{K2}
for $\overline{E_0(F/F)}$.

\subsection{Selberg case}
Let $\z_M(s)$ be the Selberg (or Ruelle) zeta function
$$
\z_M(s)=\prod_{P\in\mathrm{Prim}(M)}(1-N(P)^{-s})^{-1}
$$
of a compact Riemann surface $M$ of genus $g\ge2$, where $\mathrm{Prim}(M)$ denotes the prime geodesics on $M$
with $N(P)=\exp(\mathrm{length}(P))$. 
Let $a$, $b$, $c\in i\R$ and
$$
Z_M^{abc}(s)=\prod_{P\in\mathrm{Prim}(M)}
(1-(N(P)^a+N(P)^b)N(P)^{-s}+N(P)^{c-2s})^{-1}.
$$
Then we have the same Theorems A and B characterizing
$\z_M(s-a)\z_M(s-b)$ among $Z_M^{abc}(s)$ by using Theorem 9 of \cite[p.232]{K3}.

\subsection{More parameters}
It is possible to generalize the situation with more parameters (or representations).
For example, let $a$, $b$, $c$, $d\in i\R$ and
$$
Z^{abcd}(s)=\prod_{p:\text{ prime }}(1-(p^a+p^b+p^c)p^{-s}+
(p^{a+b}+p^{b+c}+p^{c+a})p^{-2s}-p^{d-3s})^{-1}.
$$
Then we have the following result by a similar proof:
$Z^{abcd}(s)$ is meromorphic in all $s\in\C$ iff $a+b+c=d$.
This result characterizes
$\z(s-a)\z(s-b)\z(s-c)$ among $Z^{abcd}(s)$.

Moreover, we have the following Theorem C generalizing Theorems A and B.
This characterizes $\z(s-a_1)\cdots\z(s-a_n)$ for $a_1,...,a_n\in i\R$ with $n\ge2$.
\begin{theorem}
For $n\ge2$ and $a_1,...,a_n,b\in i\R$, let
$$
Z(s)=\prod_{p:\ \text{prime}}\l((1-p^{a_1-s})\cdots(1-p^{a_n-s})+(-1)^n(p^b-p^{a_1+\cdots+a_n})p^{-ns}\r)^{-1}.
$$
Then $Z(s)$ has an analytic continuation to all $s\in\C$ as a meromorphic function if and only if
$a_1+\cdots+a_n=b$ that is $Z(s)=\z(s-a_1)\cdots\z(s-a_n)$.
When $a_1+\cdots+a_n\ne b$, it holds that $Z(s)$ is meromorphic in $\Re(s)>0$ with the natural boundary $\Re(s)=0$.
\end{theorem}

\noindent
{\it Proof.}
The method is quite similar to the case of $n=2$ treated in the proofs of Theorems A and B.

We define
$$
H(T)\in 1+TR(\R)[T]
$$
by
\begin{align*}
H(T)
&=(1-\chi_{a_1}T)\cdots(1-\chi_{a_n}T)+(-1)^n(\chi_b-\chi_{a_1+\cdots+a_n})T^n\\
&=1-(\chi_{a_1}+\cdots+\chi_{a_n})T+\cdots+(-1)^n\chi_b T^n.
\end{align*}
Then it is sufficient to show the equivalence of

(1) $a_1+\cdots+a_n=b$,

and

(2) $H(T)$ is unitary.

\bigskip
(1)$\Longrightarrow$(2):
If $a_1+\cdots+a_n=b$, then $H(T)=(1-\chi_{a_1}T)\cdots(1-\chi_{a_n}T)$, which is unitary.

(2)$\Longrightarrow$(1):
Suppose that $H(T)$ is unitary. Then, for $x\in\R$
\begin{align*}
H_x(T)
&=(1-e^{i\t_1(x)}T)\cdots(1-e^{i\t_n(x)}T)\\
&=1-(e^{i\t_1(x)}+\cdots+e^{i\t_n(x)})T+\cdots+(-1)^ne^{i(\t_1(x)+\cdots+\t_n(x))}T^n
\end{align*}
with $\t_j:\ \R\to\R$. By comparing with
$$
H_x(T)=1-(e^{a_1x}+\cdots+e^{a_nx})T+\cdots+(-1)^ne^{bx}T^n
$$
we obtain the following identities for all $x\in\R$:
\begin{align}
&e^{a_1x}+\cdots+ e^{a_nx}=e^{i\t_1(x)}+\cdots+ e^{i\t_n(x)},  \tag*{(4.\,1)}\\
&e^{(a_1+\cdots+a_{n-1})x}+\cdots+ e^{(a_2+\cdots+a_n)x}\tag*{(4.\,$n-1$)}\\
&\qquad\qquad =e^{i(\t_1(x)+\cdots+\t_{n-1}(x))}+\cdots+ e^{i(\t_2(x)+\cdots+\t_n(x))},  \nonumber\\
&e^{bx}=e^{i(\t_1(x)+\cdots+\t_n(x))},\tag*{(4.\,$n$)}
\end{align}
where $(4.\,k)$ indicates the coefficients of $T^k$ in both sides for $k=1,\,n-1,\,n$.

Now, the complex conjugate of (4.\,1) gives
\begin{align}
e^{-a_1x}+\cdots+ e^{-a_nx}=e^{-i\t_1(x)}+\cdots+ e^{-i\t_n(x)}. \tag*{($\alpha$)}
\end{align}
Dividing (4.\,$n-1$) by (4.\,$n$) we get
$$
e^{(a_1+\cdots+a_{n-1}-b)x}+\cdots+ e^{(a_2+\cdots+a_n-b)x}
=e^{-i\t_1(x)}+\cdots+ e^{-i\t_n(x)}
$$
that is
\begin{align}
e^{(a_1+\cdots+a_{n}-b)x}(e^{-a_1x}+\cdots+ e^{-a_nx})
=e^{-i\t_1(x)}+\cdots+ e^{-i\t_n(x)}. \tag*{($\beta$)}
\end{align}
Then $(\beta)-(\alpha)$ implies
$$
(e^{(a_1+\cdots+a_{n}-b)x}-1)(e^{-a_1x}+\cdots+ e^{-a_nx})=0
$$
for all $x\in\R$. Hence we obtain
$$
\f{e^{(a_1+\cdots+a_{n}-b)x}-1}x
(e^{-a_1x}+\cdots+ e^{-a_nx})=0
$$
for all $x\in\R\setminus\{0\}.$ Thus, letting $x\to0$ we have $a_1+\cdots+a_n=b$.
\hfill\qed

\end{document}